\documentclass{article}

\usepackage{PRIMEarxiv}

\usepackage[utf8]{inputenc} 
\usepackage[T1]{fontenc}    
\usepackage{hyperref}       
\usepackage{url}            
\usepackage{booktabs}       
\usepackage{amsfonts}       
\usepackage{nicefrac}       
\usepackage{microtype}      
\usepackage{lipsum}
\usepackage{fancyhdr}       
\usepackage[pdftex]{graphicx}
\graphicspath{{./pdf/}{./img/}}
\DeclareGraphicsExtensions{ .pdf, .png, .jpg, .jpeg}
\usepackage{caption}
\graphicspath{{media/}}     
\usepackage{xcolor}
\usepackage{threeparttable}
\usepackage{makecell}
\usepackage{multirow}
\usepackage{amsmath,amssymb}

\usepackage{nomencl}
\usepackage[toc,page]{appendix}
\usepackage{lineno}
\usepackage{setspace}
\usepackage{longtable}
\title{Joint Planning of Charging Stations and Power Systems for Heavy-Duty Drayage Trucks}

\author{
  Zuzhao Ye, Nanpeng Yu \\
  Department of Electrical and Computer Engineering \\
  University of California, Riverside \\
  Riverside, CA, USA\\
  \texttt{zye066@ucr.edu, nyu@ece.ucr.edu} \\
   \And
  Ran Wei \\
  School of Public Policy \\
  University of California, Riverside \\
  Riverside, CA, USA\\
  \texttt{ran.wei@ucr.edu} \\
}

\begin{document}
\doublespacing
\maketitle
\thispagestyle{fancy}

\begin{abstract}
As global concerns about climate change intensify, the transition towards zero-emission freight is becoming increasingly vital. Drayage is an important segment of the freight system, typically involving the transport of goods from seaports or intermodal terminals to nearby warehouses. This sector significantly contributes to not only greenhouse gas emissions, but also pollution in densely populated areas. This study presents a holistic optimization model designed for an efficient transition to zero-emission drayage, offering cost-effective strategies for the coordinated investment planning for power systems, charging infrastructure, and electric drayage trucks. The model is validated in the Greater Los Angeles area, where regulatory goals are among the most ambitious. Furthermore, the model's design allows for easy adaptation to other regions. By focusing on drayage trucks, this study also paves the way for future research into other freight categories, establishing a foundation for a more extensive exploration in this field.

\end{abstract}


\section{Introduction}
In the global efforts to contain climate change, a significant focus is placed on the transportation sector, especially given its contribution to greenhouse gas (GHG) emissions \cite{usdot2023climate}. In the transportation sector, medium- and heavy-duty (MDHD) vehicles, though only representing 5\% of the total vehicle fleet, are a major concern due to their contribution of over a quarter of the sector’s GHG emissions \cite{smith2020medium}. In addressing this challenge, programs such as the Advanced Clean Truck initiative, a collaborative effort spanning eighteen states in the United States led by California \cite{usdoe2023adoption}, alongside the European Union's proposal for a substantial 90\% reduction in truck emissions by 2040 \cite{mulholland2023eu}, underscore the worldwide attention to decarbonize the trucking subsector.

Drayage trucks, a specific category of heavy-duty vehicles, are primarily involved in the transport of goods between ports and intermodal terminals \cite{usepa2023env}. These trucks are predominant in areas with high population density, including disadvantaged communities, leading to not only GHG emissions but also significant air pollution and health risks \cite{carb2023fact}. California's aggressive policy to achieve full zero-emission status for drayage trucks by 2035, ahead of other heavy-duty vehicles, is a strong signal of urgency in addressing this issue \cite{newsom2020climate}. Responsive measures to this policy are already in motion: since April 2022, the ports of Long Beach and Los Angeles have started collecting \$10 per twenty-foot equivalent unit charges on non-zero-emission drayage trucks through the Clean Truck Fund initiative, aimed at spearheading the transition \cite{fisher2022socal}. Moreover, starting from 2024, only zero-emission drayage trucks will be permitted for registration in California \cite{carb2023drayage}.

Electrification is a well-recognized pathway toward zero-emission drayage \cite{miller2020electric}. However, the shift to electric drayage trucks presents substantial challenges. These vehicles, designed for long hours of operation and heavy loads, require high-capacity batteries and extensive charging infrastructure \cite{tanvir2021feasibility, kotz2022port}. The scarcity of such infrastructure, coupled with the high upfront costs of large batteries, poses significant logistical and economic challenges. This necessitates a well-thought-out approach to infrastructure development, balancing cost-effectiveness with operational efficiency. Furthermore, the power delivery system in regions with dense drayage truck traffic is already burdened by high demand from passenger, industrial, and other types of electric vehicles. Consequently, building charging infrastructure for these trucks introduces vital concerns about increased power demand, deteriorating power quality, and eventually the need for expensive power grid upgrades \cite{bradley2019charging}.

Despite various challenges, there is promising progress towards achieving zero-emission drayage. Multiple original equipment manufacturers (OEMs) are rolling out electric truck models, with a leading OEM already demonstrating a model with a daily operational range of more than 1,000 miles \cite{lamber2023tesla}, among other existing \cite{calstart2022zeti} and forthcoming models \cite{brawner2023electric}. On the charging infrastructure front, there are numerous major projects underway across the United States \cite{john2023california}. For instance, a facility equipped with 44 charging ports in Southern California \cite{adler2023forum} and a mobile charging site designed to deliver 200MW of charging power at the Port of Newark \cite{john2021what} are just two examples of the ongoing developments. The evolution of the upstream power grid is also gaining increased focus, with grid operators and utility companies actively engaging with fleet owners and charging infrastructure developers to explore a range of solutions \cite{macdougall2022as}.

The problem, however, lies in the lack of planning coordination among these sectors as they evolve at their own pace. Upgrading the power grid, which includes substations, power lines, and other critical components, is notably time-consuming, often requiring 5-15 years to perform engineering plan, obtain necessary permit, and complete construction, compared with less than two
years for new EV charging infrastructure and the fast evolution of electric truck models \cite{iea2023electricity}. The slow pace in upgrading and constructing power delivery infrastructure could significantly impede the transition to zero-emission drayage. Therefore, long-term planning of these upgrades is crucial to meet future electric drayage truck charging demands. Given the hierarchical relationship among the key elements – with charging infrastructure being the immediate downstream to the power grid, and electric vehicles being downstream to charging infrastructure – a cohesive and well-coordinated planning approach is the key \cite{unterluggauer2022electric, ahmad2022optimal}.

In this research, we aim to provide such an approach to answer the fundamental question: In line with regulatory goals \cite{carb2023fact}, what are the most cost-effective strategies for deploying charging infrastructure and upgrading power delivery systems to support the projected adoption of electric drayage trucks? To achieve this, we developed a spatio-temporal optimization model that will coordinatedly identify the most ideal locations for charging stations, and evaluate if the existing capacities of nearby electrical substations are sufficient to support this new infrastructure, or if upgrades will be required. Additionally, our model will also determine the optimal battery size for each truck and optimize charging schedules to align with time-of-use electricity pricing, thereby enhancing cost-effectiveness and truck fleet operation efficiency. In the case study, we focus on the Greater Los Angeles area, home to the Port of Long Beach and the Port of Los Angeles, two of the largest seaports globally. The proposed model utilizes a large set of real-world GPS data from over 4,800 drayage trucks, collected over a year, and pairs it with the data of more than 300 electrical substations provided by a major regional electric utility company (Southern California Edison). This approach allows for a comprehensive understanding of the specific needs and dynamics of this area, ensuring that our solution is practical and insightful.

The structure for the rest of this article is as follows: Section \ref{sec:literature} offers a review of relevant literature. In Section \ref{sec:contribution}, we highlight the key contributions of our work. Section \ref{sec:method} introduces the problem under study, detailing the formulation of the optimization framework, including its objectives, constraints, and assumptions. Section \ref{sec:data} describes the acquisition and processing of the necessary input data. The outcomes of the optimization process, along with an analysis of their impacts and implications, are presented in Section \ref{sec:results}. The article concludes with Section \ref{sec:conclusion}, summarizing our findings and proposing directions for future research.

\section{Literature Review}
\label{sec:literature}

\subsection{Charging Infrastructure Planning}

The field of charging infrastructure planning for electric vehicles is rapidly evolving, with the primary goal being to identify optimal locations and capacities for charging stations. The existing research on charging infrastructure planning can be separated into three groups based on vehicle type: passenger EVs, public transit systems, and more recently, MDHD trucks.

For passenger EVs, studies like \cite{carra2022sustainable} and \cite{roy2022examining} emphasize the importance of key high-level datasets such as demographic data, traffic flow, and EV penetration for determining optimal charging station placement and capacity. These data sets help estimate potential charging demand, which can be projected either onto geographic nodes or along travel paths. Node-based projections lead to models like the maximum covering location model, aiming to maximize demand coverage within a budget \cite{zhang2019gis, yi2022electric}, or the p-median model, focusing on minimizing cost of access to charging stations \cite{sadeghi2014optimal, janjic2021estimating}. Path-based projections focus on deploying charging stations to capture maximum EV traffic flow \cite{he2018optimal, wang2018coordinated, pal2021placement}. Unlike node- or path-based methods that estimate demand from a top-down perspective, agent-based approaches use a bottom-up strategy. They start by modeling micro-scale user behaviors, then simulate interactions between EVs, potential charging stations, and road networks to derive macro-scale outcomes, providing a more dynamic analysis of charging infrastructure needs \cite{sheppard2016cost, pagani2019user, wolbertus2021charging}. Nevertheless, agent-based methods also rely on high-level data sets such as population and EV penetration rates for fundamental simulation parameters. Generally, the planning of charging infrastructure for passenger EVs is focused on addressing the overall demand instead of the specific needs of individual vehicles.

Conversely, charging infrastructure for public transit requires detailed planning for each individual vehicle due to fixed schedules and predetermined routes. High-resolution spatio-temporal models are crucial for accurately assessing each vehicle's charging needs \cite{wei2018optimizing, stumpe2021study}. Depots and terminals are typical candidate charging station sites \cite{hsu2021depot, tzamakos2023electric}. Financial considerations usually cover both the infrastructure costs and the significant expenses associated with electric buses and batteries \cite{kunith2017electrification}, since charging infrastructure and fleets are often jointly owned by transit agencies, necessitating a comprehensive view of investment. Timetable constraints necessitate efficient charging schedules management, which can be performed by either maintaining the current fleet size \cite{ye2022decarbonizing} or optimizing the size of the fleet \cite{lee2021optimal}. Other considerations include arrival/departure uncertainties \cite{he2022integrated} and seasonal energy consumption variations \cite{liu2021optimizing}. In essence, charging infrastructure planning for public transit requires a detailed approach, addressing both charging station location/size and associated charging schedules for each individual vehicle in the fleet.

MDHD trucks have unique requirements that set them apart from passenger EVs and public transit systems. Differing from passenger EVs, trucks exhibit distinct operational patterns influenced more by local economic and spatial factors (e.g. locations of seaports, intermodal terminals, and warehouses) rather than general demographic characteristics. In contrast to public transit, truck schedules are more variable, frequently changing to accommodate the particular tasks of each day. As noted by \cite{al2021charging}, the unique operational characteristics of trucks, which include irregular timetables, longer routes, and heavier loads, present challenges not yet fully addressed in the literature. Pioneering efforts include a criteria-based method for deploying charging stations for MD and HD trucks in Minnesota \cite{khani2023identifying}. While straightforward, this method does not guarantee optimality. To mitigate this problem, a joint planning and scheduling framework specifically for drayage trucks was developed in \cite{wu2023dynamic}. However, it only determines the number of charging stations, not their geographical locations. Additionally, an agent-based approach estimating the infrastructure needs for long-haul freight in Germany was proposed in \cite{menter2023long}. The studies provide valuable insights, but the data sets used in them are relatively limited. As \cite{smith2020medium, macdougall2022as} highlighted, addressing these challenges requires extensive, large-scale data featuring broad spatial coverage, long temporal duration, substantial vehicle counts, and high frequencies.

\subsection{Joint Planning with Power System}

Parallel to infrastructure planning is the critical consideration of integrating EV charging infrastructure within existing power systems, a domain attracting increasing attention due to the surge in EV adoption. The existing literature on this topic can be split into two categories, depending on whether upgrades of the power system are considered. In scenarios where upgrades are not considered, the focus is on strategically placing charging stations in locations where existing substation capacities are adequate \cite{lin2014distribution, sadeghi2014optimal}. Additionally, operational strategies are implemented at these charging stations to minimize voltage violations and prevent transformer overloading, thus ensuring grid operation safety \cite{zhang2016pev, xiang2020reliability, deb2021novel, mao2020location}. When upgrades are factored into the planning, existing research includes the addition or reinforcement of voltage regulators, power distribution lines, transformers, and substations. These components are identified as key elements for upgrading to efficiently handle the increasing load demands \cite{arias2017robust, ehsan2019active, wang2019expansion, kabir2020demand, sa2021allocation}. As pointed out by \cite{unterluggauer2022electric}, much of the current research has been validated only on test networks. There is a recognized need to transition to large-scale, real-world case studies for more comprehensive validation. 

While researches that jointly consider power system planning was widely performed for passenger cars (e.g. \cite{davidov2019optimization, pal2023planning}) and public transit (e.g. \cite{lin2019multistage, boonraksa2019impact}), the impacts from MDHD trucks have rarely been systematically considered by the literature. In the limited studies available, \cite{bradley2019charging} introduces a set of algorithms for the strategic placement of charging stations for drayage trucks, considering grid constraints for short-term planning and assuming grid capacity sufficiency in the long term. Meanwhile, \cite{londono2019optimal} addresses energy losses in the planning of charging infrastructure for freight vehicles, and \cite{khani2023identifying} identifies substation proximity as a key factor in choosing locations for truck charging stations. A particular point of interest is the integration of heavy-duty drayage trucks within these systems. Their distinct high-power charging needs, particularly in densely populated urban regions, present substantial challenges to the capacities of existing power delivery systems. This situation underscores an immediate need for dedicated research and the development of innovative infrastructure solutions, areas that are currently underrepresented in existing literature.

\subsection{Contribution}
\label{sec:contribution}

Reflecting on the current status of research, we underscore the significance of this study's contributions as follows:

\begin{itemize}
\item While extensive research has been conducted on charging infrastructure planning for passenger EVs and public transit, similar studies for MDHD trucks are lacking, particularly for HD drayage trucks, which are major contributors to urban pollution. This study aims to bridge this gap, focusing on HD drayage trucks as an example subset of MDHD trucks, thereby shedding light on planning approaches for other types of trucks.

\item Considering the substantial charging power needed by HD drayage trucks in densely populated areas, our study explicitly considers the load hosting capacity of the nearby power delivery system. We have developed an integrated framework that jointly determines the placement of charging stations, their grid connections, and any necessary substation upgrades.

\item Utilizing three large-scale realistic data sets – truck GPS trajectories, potential charging locations, and existing electric substations – we conduct a detailed case study in the Greater Los Angeles area, home to two major seaports. This analysis evaluates the adequacy of the current power system capacity and offers practical, cost-effective strategies for meeting state regulatory goals.

\end{itemize}

\section{Methodology}
\label{sec:method}

In this section, we introduce the problem addressed in this study, detailing the specific modeling and constraints relevant to each sector: Drayage trucks, charging infrastructure, and the power delivery system. Our model is designed to be solved under two distinct objective modes: Examining hosting capacity or ensuring compliance with regulatory goals in the most cost-effective manner.

\begin{figure*}[!h]
\centering
\captionsetup{justification=centering}
\includegraphics[width=\textwidth, clip=true, trim = 0mm 0mm 0mm 0mm]{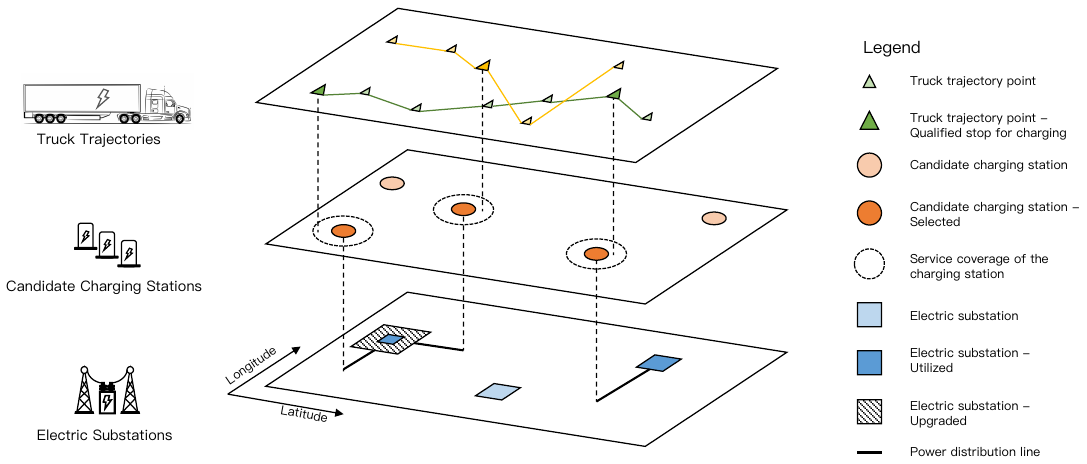}
\caption{An illustration of the systems studied.}
\label{fig:overview}
\end{figure*}

\subsection{Problem Description}
\label{sec:problem-description}

We consider a regional area populated by a number of fossil fuel-powered drayage trucks, denoted as the set $I$. To support their transition to electric vehicles, a set $J$ of potential charging stations, including depots, intermodal terminals, and existing truck stops, has been identified. These charging stations will rely on the nearby power delivery systems, represented by a set $K$ of electric substations. The drayage trucks will make use of strategically selected charging stations along their routes to ensure sustained operations during the time horizon $T$, with these charging stations sourcing their power from the power delivery system, as illustrated in Figure \ref{fig:overview}. This study is based on the following assumptions:

\begin{itemize}
\item The routes, or trajectories, of the trucks remain unchanged before and after their transition to electric vehicles.
\item A truck can use a charging station if the station is selected for deployment and is either its designated depot or a public charging facility.
\item A truck is eligible for charging if its trajectory indicates a parking duration of at least 30 minutes.
\end{itemize}

Different stakeholders in the system have their unique concerns and interests. For example:

\begin{itemize}
\item Fleet owners are focused on whether delivery tasks can be efficiently completed, strategies for managing fleet transitions, and optimal scheduling for charging.
\item Charging infrastructure developers prioritize decisions about the placement of charging stations and determining their appropriate scale.
\item Electric utilities are primarily concerned with evaluating if the existing load hosting capacity of substations can handle the increased demand and identifying where upgrades might be required.
\item Regulators are interested in gaining a holistic view of the entire system to develop informed policies and to strategically fund essential components of the transition.
\end{itemize}
Each group also shares a common concern about the potential costs involved. In the following sections, we will introduce the constraints specific to each sector and outline objectives that address these varied concerns. Table \ref{tab:params} summarizes the notations used in this study.

\begin{longtable}{ll}
\caption{Summary of notations for sets, parameters, and decision variables.} \label{tab:params} \\
\toprule
 & \textbf{Sets} \\
\midrule
\endfirsthead

\multicolumn{2}{c}%
{\tablename\ \thetable\ -- \textit{Continued from previous page}} \\
\toprule
 & \textbf{Parameters} \\
\midrule
\endhead

\midrule
\multicolumn{2}{r}{\textit{Continued on next page}} \\
\endfoot

\bottomrule
\endlastfoot

$I$ & Set of drayage trucks.\\
$J$ & Set of candidate charging stations.\\
$J_{it}$ & Subset of $J$, accessible charging stations for truck $i$ during time period $t$.\\
$J_k$ & Subset of $J$, charging stations covered by electric substation $k$.\\
$K$ & Set of electric substations.\\
$K_j$ & Subset of $K$, substations accessible to charging station $j$.\\
$T$ & Set of discrete time periods in the study time horizon.\\

\midrule
 & \textbf{Parameters} \\
\midrule
$a_{ijt}$ & Binary parameter: 1 if truck $i$ can access charging station $j$ during time period $t$, else 0.\\
$c^{veh}$ & Amortized cost of a base electric truck model.\\
$c^{btr}$ & Amortized cost for each extra kWh of battery capacity.\\
$c^{ele}_t$ & Cost per kWh of electricity at time $t$.\\
$c^{ctr}$ & Amortized cost of charging station construction, including planning, permit, civil work, etc.\\
$c^{cap}$ & Amortized cost of establishing power capacity at a charging station per kW basis.\\
$c^{chg}$ & Amortized cost of charger per kW basis.\\
$c^{lne}$ & Amortized cost per mile for a distribution line between charging station $j$ and substation $k$.\\
$c^{upg,std}$ & Amortized standard fixed cost of substation upgrade.\\
$c^{upg,var}$ & Amortized variable cost for additional substation capacity upgrade.\\
$D$ & The number of days in a year.\\
$d_{jk}$ & Distance between charging station $j$ and electric substation $k$.\\
$\Delta e^{con}_{it}$ & Energy consumption of truck $i$ at time period $t$.\\
$e^{bse}$ & Base battery capacity of an electric truck.\\
$g_t$ & Grid electricity carbon intensity at time $t$.\\
$\kappa$ & Battery charging cycle efficiency.\\
$N_y$ & Target number of electric drayage trucks by year $y$.\\
$p^{max}$ & Maximal charging power for a truck.\\
$P^{sbs}_{k}$ & Remaining load hosting capacity of electric substation $k$.\\
$P^{upg,std}$ & Standard capacity increment for a substation upgrade.\\
$PF$ & Typical power factor.\\
$r^{stp}_{it}$ & The proportion of time period $t$ that the truck $i$ is stopped. \\
$u_{i}$ & Daily carbon emission of truck $i$ if it uses diesel as the fuel.\\

\midrule
 & \textbf{Decision Variables$^*$} \\
\midrule
$\gamma_{jk}$ & (B) Indicates if station $j$ is connected to substation $k$.\\
$e_{it}$ & (C) Battery energy level of truck $i$ at time $t$.\\
$e^{cap}_{i}$ & (C) Required battery capacity of truck $i$.\\
$p_{it}$ & (C) Charging power for truck $i$ at time $t$.\\
$p^{trk}_{ijt}$ & (C) Charging power supplied by charging station $j$ to truck $i$ at time $t$.\\
$\hat p^{trk}_{ijt}$ & (B) Indicates if truck $i$ charges at charging station $j$ at time $t$.\\
$P_{jt}$ & (C) Total power demand of trucks in charging station $j$ at time $t$.\\
$P^{chs}_{j}$ & (C) Power capacity of charging station $j$.\\
$\hat P^{chs}_{j}$ & (B) Indicates if candidate charging station $j$ is deployed.\\
$P^{chs-sbs}_{jk}$ & (C) Peak power flow between charging station $j$ and electric substation $k$.\\
$P^{upg}_{k}$ & (C) Total capacity upgrade at substation $k$.\\
$P^{upg,var}_{k}$ & (C) Variable part of capacity upgrade at substation $k$.\\
$\hat P^{upg}_{k}$ & (B) Indicates an upgrade at substation $k$.\\
$x_{i}$ & (B) Indicates if truck $i$ is converted to electric. 
\\

\end{longtable}
\vspace{-0.7cm}
\noindent
\textsuperscript{*}\textit{(C) Continuous variable; (B) Binary variable.}

\subsection{Drayage Trucks}
In this subsection, we introduce two categories of constraints relevant to the drayage truck sector. First, we examine the battery energy constraints applicable to trucks if they are electrified. Second, we discuss the spatial-temporal constraints that govern their charging activities. Additionally, we will discuss the detailed costs for this sector.

\textbf{Battery Energy Constraints}: Constraint \eqref{eq:cons-truck-1} defines the transition rule of the battery energy level $e$ of the truck $i$ from time $t$ to the next time period $t^{\prime}$, where $p_{it}$ and $\Delta e^{con}_{it}$ is the charging power (in kW) and the energy consumption (in kWh) during this time period, respectively. $r^{stp}_{it}$ is the proportion of time that the truck is in a stopping condition during time period $t$. Charging-related energy loss is accounted for by $(1-\sqrt{\kappa}) \cdot p_{it}$, where $\kappa$ represents battery cycle efficiency \cite{foggo2017improved}. The $Next(\cdot)$ function operates as follows: if $t$ is the final period in the time horizon $T$, $Next(t)$ reverts to the first period of $T$; otherwise, $Next(t) = t+\Delta t$. This creates a cyclical energy loop for the truck, ensuring continuous operation. Constraint \eqref{eq:cons-truck-2} ensures the battery's state-of-charge (SoC) remains within a specific range, where $e^{cap}_i$ is truck $i$’s battery capacity. Constraint \eqref{eq:cons-truck-3} introduces binary variables $x_{i}$ to indicate whether the truck $i$ is converted from fossil fuel-powered to electric.

\begin{equation}
    e_{it^{\prime}} = e_{it} + \left (1-\sqrt{\kappa} \right ) \cdot p_{it} \cdot r^{stp}_{it} \cdot \Delta t - x_i \cdot \Delta e^{con}_{it} , \quad \forall i \in I, \forall t \in T, t^{\prime} = Next(t)
\label{eq:cons-truck-1}
\end{equation}

\begin{equation}
SoC^{min} \cdot e^{cap}_i \leq e_{it} \leq SoC^{max} \cdot e^{cap}_i, \quad \forall i \in I, \forall t \in T,
\label{eq:cons-truck-2}
\end{equation}

\begin{equation}
    x_i \in \{0, 1\}, \quad \forall i \in I,
\label{eq:cons-truck-3}
\end{equation}

\textbf{Charging Activity Spatial-Temporal Constraints}: The charging power for truck $i$ at time $t$ at a candidate charging station $j$ is represented as $p^{trk}_{ijt}$. Constraint \eqref{eq:cons-truck-4} sets the boundaries for $p^{trk}_{ijt}$, where its lower limit is zero, and its upper limit depends on two factors: if truck $i$ is an electric truck, and if it has access to charging station $j$ during time $t$, as indicated by the parameter $a_{ijt}$. In our study, a truck is considered to have access to station $j$ if it is within 0.5 mile of the station at time $t$, thereby setting $a_{ijt}=1$ in such instances.

\begin{equation}
    0 \leq p^{trk}_{ijt} \leq x_i \cdot a_{ijt} \cdot p^{max}, \quad \forall i \in I, \forall t \in T, \forall j \in J
\label{eq:cons-truck-4}
\end{equation}

The charging power of truck $i$ at time $t$, denoted as $p_{it}$, is calculated as the sum of $p^{trk}_{ijt}$ across all relevant stations, as defined in constraint \eqref{eq:cons-truck-7}.

\begin{equation}
    p_{it} = \sum_{j \in J_{it}}p^{trk}_{ijt}, \quad \forall i \in I, \forall t \in T,
\label{eq:cons-truck-7}
\end{equation}

where $J_{it} = \{j | j \in J \text{ and } a_{ijt}=1\}$ is a subset of $J$ consisting of the charging stations accessible to truck $i$ during the time period $t$.

Additionally, each truck $i$ is restricted to using at most one charging station at any given time. To implement this logical constraint, we introduce a binary variable $\hat p^{trk}_{ijt}$, as detailed in \eqref{eq:cons-truck-5}, which is set to 1 whenever there is non-zero charging power. Constraint \eqref{eq:cons-truck-6} then ensures that a truck utilizes no more than one charging station simultaneously.

\begin{equation}
\hat p^{trk}_{ijt} =
\begin{cases}
    1, & \text{if $p^{trk}_{ijt} > 0$}\\
    0, & \text{otherwise},
\end{cases}, \forall (i, j, t) \in (I, J_{it}, T)
\label{eq:cons-truck-5}
\end{equation}

\begin{equation}
\sum_{j \in J_{it}} \hat p^{trk}_{ijt} \leq 1, \quad \forall i \in I, \forall t \in T.
\label{eq:cons-truck-6}
\end{equation}

When truck $i$ is parked and charging at station $j$ during time $t$, it is practically expected to remain at this station for the entire charging session, implying that charging should occur at a single location continuously. Constraint \eqref{eq:cons-truck-8} enforces this logic. The first term on the left-hand side of \eqref{eq:cons-truck-8} signifies a charging activity at station $j$ during time $t$, while the second term represents charging at different locations during the subsequent period $t^{\prime}$.

\begin{equation}
\hat p_{ijt} + \sum_{j^{\prime} \in J_{it^{\prime}}\setminus{j}} \hat p_{ij^{\prime}t^{\prime}} \leq 1, \quad \forall (i, j, t) \in (I, J_{it} \cap J_{it^{\prime}}, T), t^{\prime} = Next(t)
\label{eq:cons-truck-8}
\end{equation}

\textbf{Cost of the Drayage Truck Sector}: The costs in this sector, as detailed in \eqref{eq:cons-truck-9}, are split into vehicle and operational expenses. For vehicle costs, we acknowledge the diverse battery capacity needs of trucks. Those covering shorter distances or with regular access to charging may require smaller batteries, and this variability is reflected in the cost structure. We consider a fixed cost $c^{veh}$ for a truck with base battery capacity $e^{bse}$ and a variable cost $c^{btr} \cdot \left(e^{cap}_{i}-e^{bse}\right)$ based on the actual battery capacity chosen, allowing for customization to each truck’s operational requirements.

In terms of operational costs, our primary focus is on electricity expenses, which are a key differentiator between the costs of fossil-fueled and electric trucks. As outlined in \eqref{eq:cons-truck-9}, $D$ is the number of days in a year, and $c^{ele}_t$ refers to the time-of-use electricity price at time $t$ (as detailed in Section \ref{sec:data-cost}). Typically, in fossil fuel-powered trucks, fuel costs account for about 28\% of the total operational expenses, a close second to driver wages (32\%) \cite{leslie2023analysis}. Although the energy costs for electric trucks are already much lower than those for conventional trucks \cite{carb2018battery}, further reductions in electricity expenses can be achieved by optimizing charging times.

\begin{equation}
    C^{trk} = \underbrace{\sum_{i \in I}\left[c^{veh} \cdot x_i + c^{btr} \cdot \left(e^{cap}_{i}-e^{bse}\cdot x_i\right)\right]}_{\text{Cost of electric truck and battery upgrades}}
    +
    \underbrace{D\sum_{t \in T} \sum_{i \in I}c^{ele}_t \cdot p_{it}\Delta t}_{\text{Cost of electricity}}.
\label{eq:cons-truck-9}
\end{equation}

\subsection{Charging Infrastructure}
This subsection covers the power capacity needs, investment decisions, and costs associated with the candidate charging stations.

\textbf{Power Capacity Constraints}: The total power consumption $P_{jt}$ at a charging station $j$ at time $t$ is the cumulative charging load from trucks with access, as defined in \eqref{eq:cons-CS-1}. The subset $I_{jt}$ of $I$, where $I_{jt}={i| i\in I \text{ and } a_{ijt}=1}$, represents the trucks at station $j$ at time $t$. The total load at any time should not exceed the station's power capacity $P^{chs}_{j}$, as stated in \eqref{eq:cons-CS-2}.

\begin{equation}
    P_{jt} = \sum_{i \in I_{jt}} p^{trk}_{ijt}, \quad \forall j \in J, \forall t \in T,
\label{eq:cons-CS-1}
\end{equation}

\begin{equation}
    P_{jt} \leq P^{chs}_{j} , \quad \forall j \in J, \forall t \in T,
\label{eq:cons-CS-2}
\end{equation}

\textbf{Investment Decision Constraints}: A binary variable $\hat P^{chs}_{j}$ indicates the selection of a candidate charging station. Per \eqref{eq:cons-CS-3}, where $M$ is a large number, $\hat P^{chs}_{j}=1$ implies the station is chosen if its power capacity is above zero.

\begin{equation}
\hat P^{chs}_{j} \geq P^{chs}_{j}/M, \quad \forall j \in J
\label{eq:cons-CS-3}
\end{equation}

\begin{equation}
\hat P^{chs}_{j} \in \{0, 1\}, \quad \forall j \in J
\label{eq:cons-CS-4}
\end{equation}

\textbf{Cost of Charging Infrastructure}: The cost of charging infrastructure includes three parts: fixed construction costs (including planning, permitting, and civil work), expenses for establishing power capacity (such as transformers and conduits), and the cost of chargers. We approximate the costs of chargers to be proportionate to the station's power capacity, with estimates interpolated from the stepwise power specifications of chargers as outlined in \cite{crc2023access}.

\begin{equation}
    C^{chg} = \sum_{j \in J} \left [c^{ctr} \cdot \hat P^{chs}_{j} + (c^{cap}+ c^{chg}) \cdot P^{chs}_{j} \right]
\label{eq:cons-CS-5}
\end{equation}

\subsection{Power Delivery System}
This subsection explores the essential connections between charging stations and the power delivery system, along with the required upgrades for electric substations. 

\textbf{Electric Substation - Charging Station Connection Constraints}: To source electricity from the grid, each charging station must be connected to the power delivery system. Given the high-power demands of truck charging stations, each charging station is required to establish a new connection with an electric substation rather than tapping into an existing feeder. We denote the peak power draw from a substation $k$ by a charging station $j$ as $P^{chs-sbs}_{jk}$, and $K_j$ is the subset of $K$ indicating the substations accessible to charging station $j$ (access criteria detailed in Section \ref{sec:data-substations}). Constraint \eqref{eq:cons-sub-1} ensures that a charging station's power capacity is satisfied from the substations it can access. A binary variable $\gamma_{jk}$ is used to indicate a connection between charging station $j$ and substation $k$. For practical purposes, each charging station is limited to a single substation connection, as enforced by \eqref{eq:cons-sub-3}.

\begin{equation}
    P^{chs}_{j} \leq \sum_{k \in K_j} P^{chs-sbs}_{jk}, \quad \forall j \in J,
\label{eq:cons-sub-1}
\end{equation}

\begin{equation}
\gamma_{jk} =
\begin{cases}
    1, & \text{if $P^{chs-sbs}_{jk} > 0$}\\
    0, & \text{otherwise},
\end{cases}, \forall j \in J, \forall k \in K_j
\label{eq:cons-sub-2}
\end{equation}

\begin{equation}
    \sum_{k \in K_j} \gamma_{jk} \leq 1, \quad \forall j \in J,
\label{eq:cons-sub-3}
\end{equation}

\textbf{Electric Substation Hosting Capacity Constraints}: Each electric substation is limited by its remaining load hosting capacity $P^{sbs}_{k}$. The combined power demand from all connected charging stations must not exceed this limit, plus any additional capacity from potential upgrades, as defined in \eqref{eq:cons-sub-4}. Practically, any upgrade $P^{upg}_{k}$ involves installing a transformer with a standard MVA rating $P^{upg,std}$, representing the baseline upgrade capacity. Hence, $P^{upg}_{k}$ consists of this fixed standard upgrade plus a variable component, as outlined in \eqref{eq:cons-sub-5}. Here, $PF$ denotes the typical power factor, and $\hat P^{upg}_{k}$ is a binary indicator of whether an upgrade is initiated, with a value of 1 signifying an upgrade. It is also important to note that the variable component of an upgrade is only relevant when $\hat P^{upg}_{k}=1$, as stated in \eqref{eq:cons-sub-7}, with $M$ being a suitably large number.

\begin{equation}
    \sum_{j \in J_k} P^{chs-sbs}_{jk} \leq P^{sbs}_{k} + P^{upg}_{k}, \quad \forall k \in K,
\label{eq:cons-sub-4}
\end{equation}

\begin{equation}
    P^{upg}_{k} = P^{upg,std} \cdot PF \cdot \hat P^{upg}_{k} + P^{upg,var}_{k}, \quad \forall k \in K,
\label{eq:cons-sub-5}
\end{equation}

\begin{equation}
    \hat P^{upg}_{k} \in \{0, 1\}, \quad\forall k \in K,
\label{eq:cons-sub-6}
\end{equation}

\begin{equation}
    0 \leq P^{upg,var}_{k} \leq M \cdot \hat P^{upg}_{k}, \quad\forall k \in K.
\label{eq:cons-sub-7}
\end{equation}

\textbf{Cost of Power Delivery System}: The cost for the power delivery system consists of two primary components. The first is the cost of power distribution lines connecting charging stations to substations, calculated based on $d_{jk}$, the distance between them, and a per-mile cost $c^{lne}$. The second component involves the cost of upgrading substations. This includes a standard fixed cost $c^{upg,std}$ incurred for each upgrade. Additionally, if extra load hosting capacity beyond the standard upgrade is needed, the cost will be proportional to this additional capacity, with $c^{upg,var}$ serving as the unit cost per MW.

\begin{equation}
    C^{pwr} = \underbrace{\sum_{j \in J}\sum_{k \in K_j}c^{lne} \cdot \gamma_{jk} \cdot d_{jk}}_{\text{Power distribution line}}
    + \ \underbrace{\sum_{k \in K}\left(c^{upg,std} \cdot \hat P^{upg}_{k} + c^{upg,var} \cdot P^{upg,var}_{k}\right)}_{\text{Substation upgrade}}
\label{eq:cons-sub-8}
\end{equation}

\subsection{Objective Functions}
Addressing different needs of diverse stakeholders in the electric ecosystem, this optimization model can be paired with various objective functions. Here, we focus on two specific objectives that assist utility companies in evaluating current hosting capacity and guide regulators in identifying cost-effective pathways to achieve zero-emission drayage. These objectives are categorized into two distinct modes.

\subsubsection{Mode 1 - Hosting Capacity Examination}
\label{sec:method-obj-mode1}

This mode aims to maximize the number of electric drayage trucks, excluding considerations for substation upgrades:

\begin{equation}
    \underset{}{\text{max}} \ \sum_{i \in I}x_i,
\label{eq:obj-1}
\end{equation}

\begin{equation}
    P^{upg}_{k} = 0, \quad\forall k \in K
\label{eq:cons-sub-add}
\end{equation}

The optimization model is then formulated with objective function \eqref{eq:obj-1} and under constraints \eqref{eq:cons-truck-1}-\eqref{eq:cons-truck-8}, \eqref{eq:cons-CS-1}-\eqref{eq:cons-CS-4}, \eqref{eq:cons-sub-1}-\eqref{eq:cons-sub-4}, and \eqref{eq:cons-sub-add}. This mode helps utility companies assess the potential for electrifying trucks under existing grid capacities, without upgrades. The outcome offers a comparison against regulatory goals, highlighting the potential need for infrastructure improvements.

\subsubsection{Mode 2 - Cost-Effective Regulation Compliance}
\label{sec:method-obj-mode2}

Here, the goal is to minimize the cumulative costs across all sectors while meeting yearly regulatory goals:

\begin{equation}
    \underset{}{\text{min}} \ C^{trk} + C^{chg} + C^{pwr},
\label{eq:obj-2}
\end{equation}

\begin{equation}
    \sum_{i \in I}x_i \geq N_y,
\label{eq:cons-truck-add}
\end{equation}

where $N_y$ represents the desired number of electric drayage trucks by year $y$. The optimization model is formulated with the objective function \eqref{eq:obj-2} and under constraints \eqref{eq:cons-truck-1}-\eqref{eq:cons-sub-8}, and \eqref{eq:cons-truck-add}. This mode assists regulators in determining the most efficient pathways to zero-emission drayage. Additionally, GHG emissions post-conversion are calculated using \eqref{eq:cons-emission}, where $U^{trk}_y$ is the total emission from trucks in year $y$, $u^{trk}_i$ is the emission from truck $i$ if it is fossil-fueled, and $\sum_{t \in T}g_t \cdot p_{it} \cdot \Delta t$ estimates the emissions from electric charging based on the electricity carbon intensity $g_t$.

\begin{equation}
    U^{trk}_y = \sum_{i \in I}(1-x_i)\cdot u^{trk}_i + \sum_{i \in I}\sum_{t \in T}g_t \cdot p_{it} \cdot \Delta t,
\label{eq:cons-emission}
\end{equation}

\section{Data Description}
\label{sec:data}

This section introduces the data utilized in this study, which primarily includes drayage truck trajectories, candidate charging station locations, and regional electric substations. Additionally, we detail the cost estimates associated with each sector. Our research concentrates on the Greater Los Angeles area, home to two of the world's busiest seaports: the Ports of Long Beach and Los Angeles.

\subsection{Drayage Trucks}

\begin{figure}[!h]
\centering
\captionsetup{justification=centering}
\includegraphics[width=0.90\textwidth, clip=true, trim = 0mm 0mm 0mm 0mm]{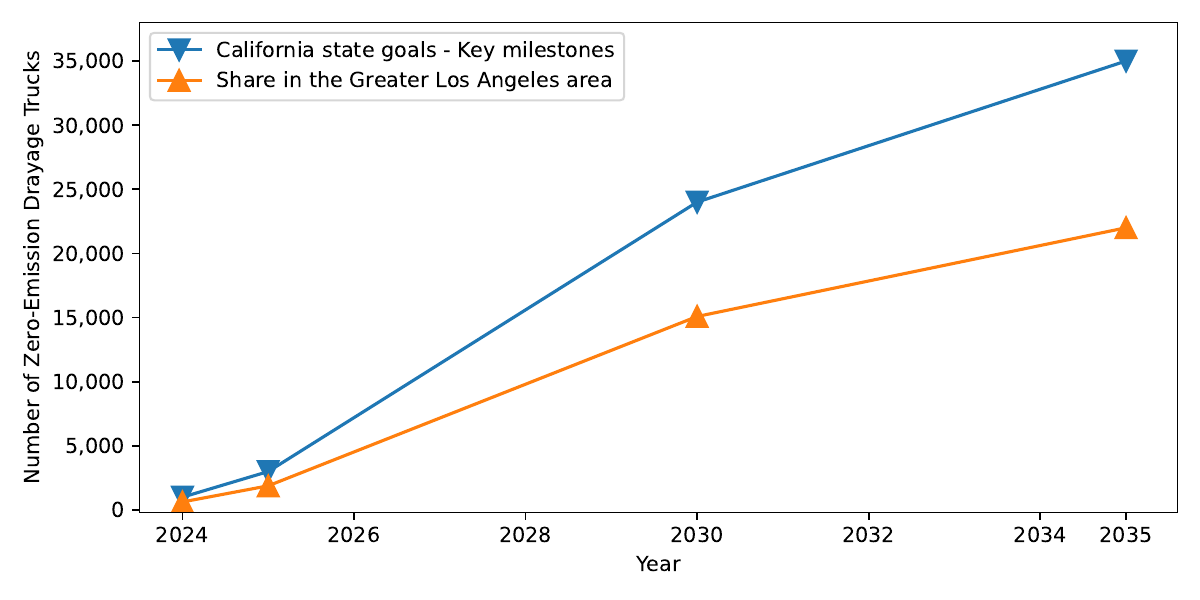}
\caption{Trends of zero-emission drayage trucks: California state goals and the share in the Greater Los Angeles area.}
\label{fig:goals}
\end{figure}

California's regulations mandate that by 2035, all drayage trucks, totaling approximately 35,000, must be zero-emission. The California Air Resources Board (CARB) has set several key milestones for this transition: 1,000 trucks by 2024, 3,000 by 2025, a significant increase to 24,000 by 2030, and ultimately 35,000 by 2035 \cite{carb2023fact}. In the Greater Los Angeles area, which includes about 22,000 of the state's drayage trucks \cite{polb2023clean}, the transition to zero-emission vehicles is expected to mirror the state's trend. Based on this assumption, we have projected the growth of zero-emission trucks in Greater Los Angeles, as illustrated in Figure \ref{fig:goals}. In our analysis, we mainly consider that most zero-emission trucks will be battery electric vehicles, given the current industry trends towards improved efficiency, lower battery costs, and expanding charging infrastructure.

Our dataset comprises detailed GPS trajectories for over 4,800 drayage trucks, spanning more than a year. The data is de-identified and includes only geocoordinates and timestamps. We transformed this raw GPS data into 15-minute frequency activities following a three-step procedure outlined in Appendix \ref{ap:gps}.

\begin{figure}[!h]
\centering
\captionsetup{justification=centering}
\includegraphics[width=1.0\textwidth, clip=true, trim = 0mm 0mm 0mm 0mm]{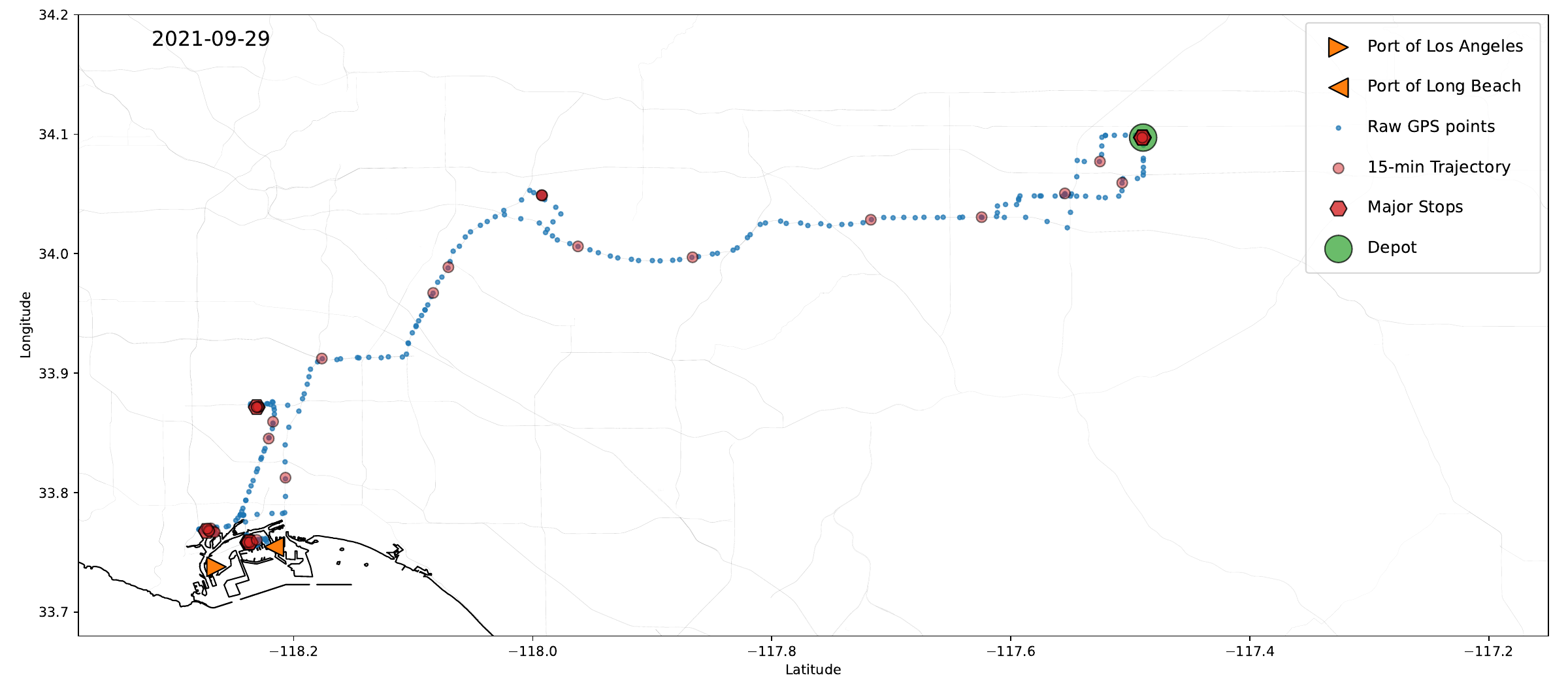}
\caption{Comparison between raw GPS data and processed 15-minute trajectories (Depot location adjusted for confidentiality).}
\label{fig:trajectories}
\end{figure}

\begin{figure}[!h]
\centering
\captionsetup{justification=centering}
\includegraphics[width=1.0\textwidth, clip=true, trim = 0mm 0mm 0mm 0mm]{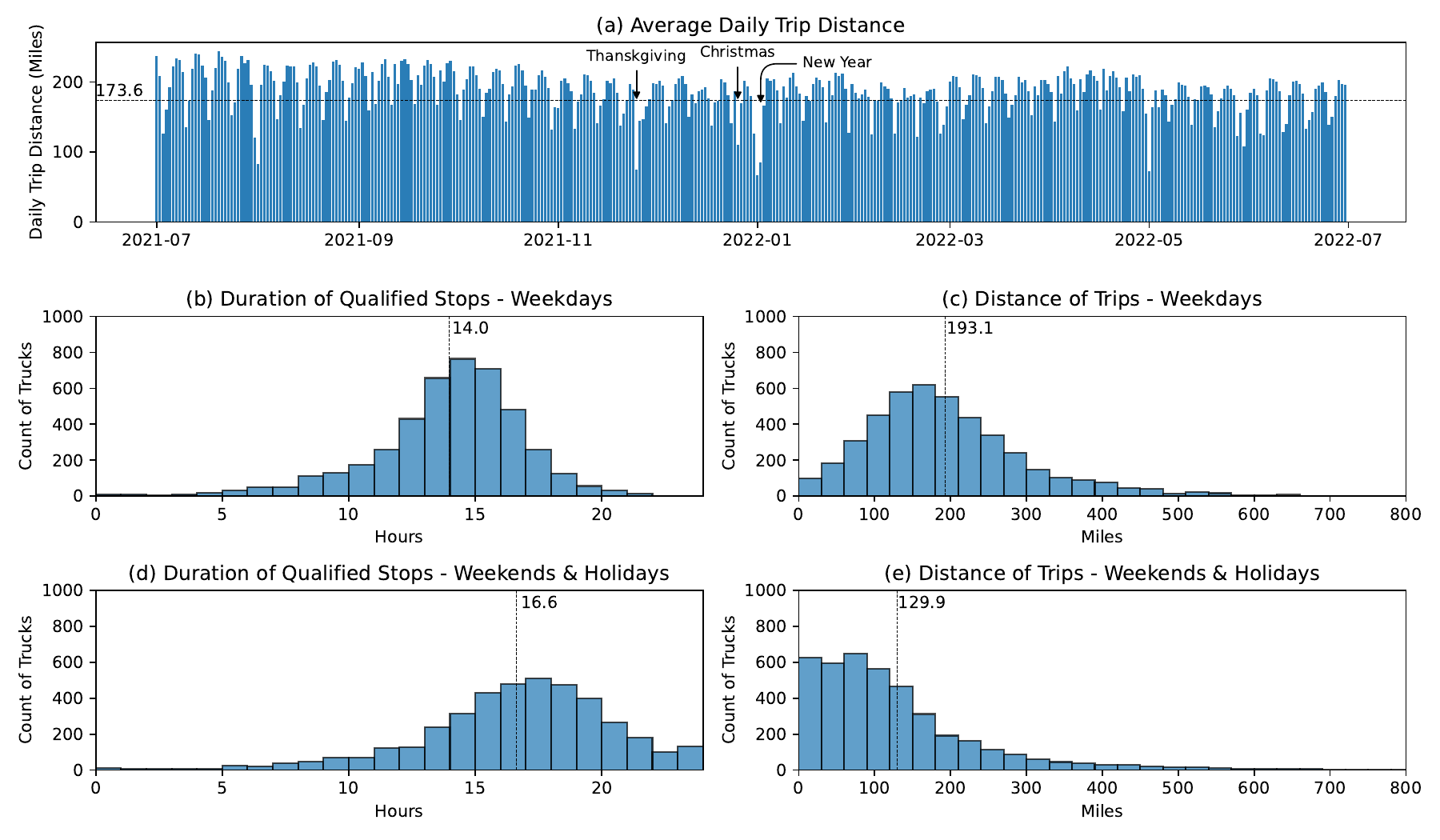}
\caption{Truck activity statistics with dashed lines indicating average values for each category.}
\label{fig:stats}
\end{figure}

During the data processing phase, we also compiled key statistics. Figure \ref{fig:stats}(a) displays the average daily travel distances for the trucks in our dataset, calculated daily throughout the data collection period. This data shows a distinct weekly pattern, with noticeable reductions in travel distances on weekends and holidays. Additionally, we created histograms, as shown in Figures \ref{fig:stats}(b)-(e), which detail the duration of qualified stops (defined as continuous stops of 30 minutes or longer) and daily travel distances for weekdays and weekends, calculated on an individual truck basis. These histograms indicate that, on average, trucks in our dataset make stops of about 14.0 hours per day on weekdays and 16.6 hours on weekends, counting only stops longer than 30 minutes. Meanwhile, the average daily travel distance is 193.1 miles on weekdays, which reduces to 129.9 miles on weekends. These statistics are in good alignment with the findings reported in \cite{carb2020advanced, kotz2022port}.

Our objective is to analyze a typical day, for which we selected September 29, 2021—a weekday with relatively high daily travel distances. To refine our dataset, we applied several filters: including only trucks with data on this date, limiting to those operating within the Greater Los Angeles area, and excluding trucks traveling less than 10 miles to omit inactive ones. This filtering resulted in a dataset of 733 trucks, which we believe are representative of local drayage activities. To approximate the actual fleet size of 22,000 drayage trucks in Greater Los Angeles, we replicated these 733 truck trajectories 30 times, creating a comprehensive set $I$ of 21,990 trucks. It is important to note that while the trajectories were duplicated, their energy consumptions were varied by multiplying with a random factor uniformly drawn from [0.95, 1.05], to account for variations in driving behavior.

\begin{figure}[!h]
\centering
\captionsetup{justification=centering}
\includegraphics[width=1.0\textwidth, clip=true, trim = 0mm 0mm 0mm 0mm]{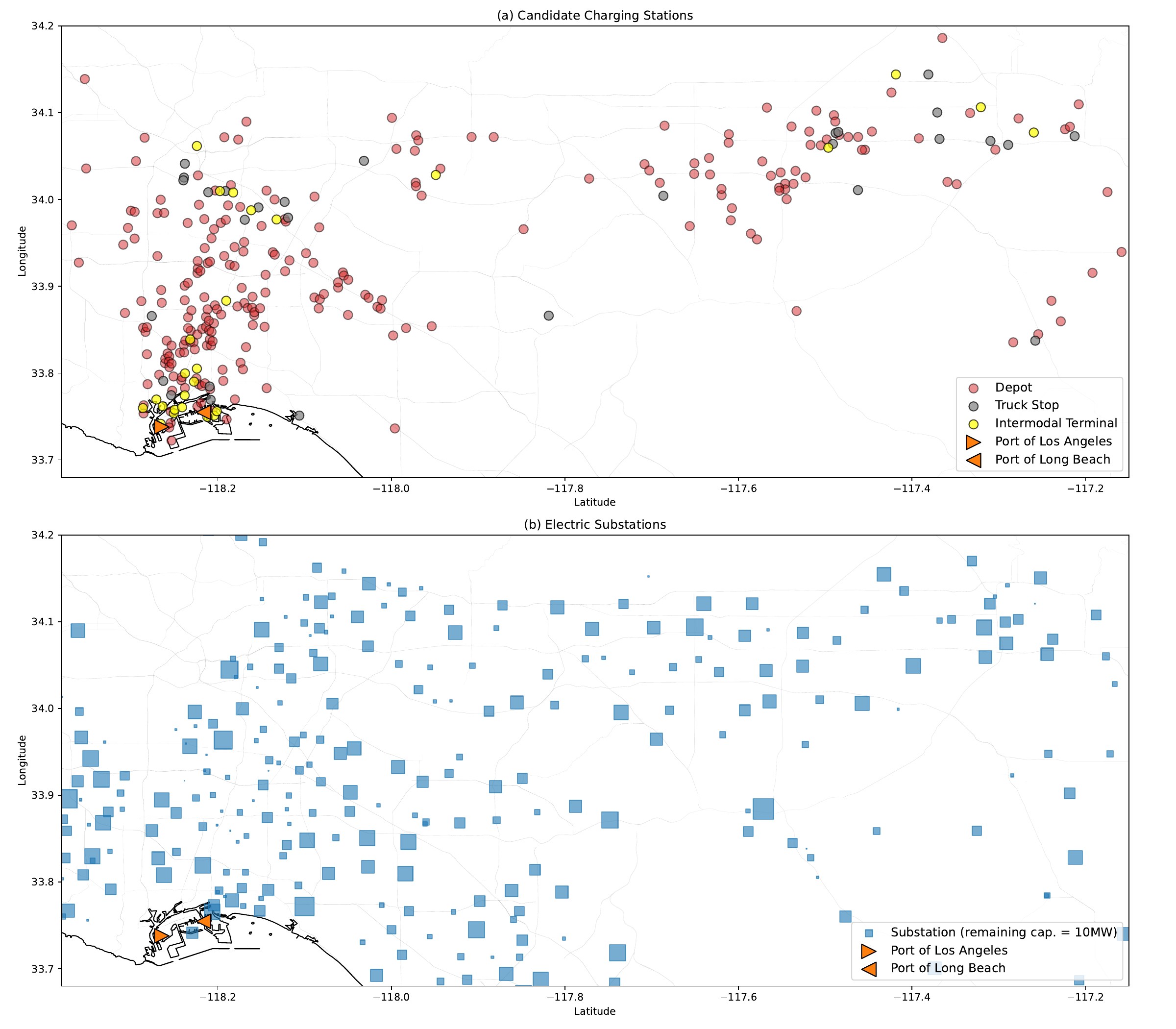}
\caption{Geographical distribution of candidate charging stations and electric substations.}
\label{fig:CS-subs-combined}
\end{figure}

\subsection{Candidate Charging Stations}

This study identifies three types of locations as potential charging stations: depots, truck stops, and intermodal terminals. Depots are determined from truck trajectories within the month of the representative day, with the longest stop presumed to be the depot. Using the DBSCAN clustering algorithm \cite{ester1996density}, known for efficiently grouping points based on density, we merge nearby depots within a 1,000 feet (304.8 meters) radius. Existing truck stop locations are sourced from \cite{allstays2023truck}, and intermodal terminals from \cite{loadmatch2023intermodal}. This process results in 262 depots, 55 truck stops, and 27 intermodal terminals, forming a total of 344 candidate charging stations (set $J$). Figure \ref{fig:CS-subs-combined}(a) illustrates their geographical spread, showing a concentration around the port area with extensions further inland, aligning with the distribution of local warehouses \cite{warehouse2023}. To determine truck access to these stations, we overlap the 15-minute trajectories with station locations. A truck is considered to have access if its qualified stop is within 0.5 miles of a public station (truck stop or intermodal terminal) or its depot, and with this assumption, we derived the list of parameters $a_{ijt}$. Given the rapid advancements in charging technology, including the introduction of megawatt chargers \cite{crc2023access}, we set $p^{max}$, the maximum charging power for a truck, at 1MW to align with these evolving standards.

\subsection{Electric Substations}
\label{sec:data-substations}
Southern California Edison (SCE) predominantly powers the Greater Los Angeles area. We utilize publicly available data on SCE's electrical distribution substations, including estimated remaining load hosting capacities \cite{drpep}. Our study includes 255 substations within the study area. Figure \ref{fig:CS-subs-combined}(b) shows these substations, with circle sizes representing their remaining load hosting capacities. Unlike the charging stations, substations are more uniformly distributed to serve various load types. We assume that each candidate charging station can access the five nearest substations. It should be noted that while this study does not include the construction of new substations, such additions could be easily accommodated once candidate locations are provided.



\subsection{Cost Estimates}
\label{sec:data-cost}

\begin{table}[ht]
\doublespacing
\caption{List of Initial Investments, Lifespans, and the Amortized Annual Costs}
\begin{center}
\begin{threeparttable}
\begin{tabular}{p{0.25\textwidth} | p{0.2\textwidth} | p{0.08\textwidth} | p{0.27\textwidth}}
    \toprule
    \textbf{Asset} & \textbf{Initial Investment} & \textbf{Lifespan} & \textbf{Annual Cost} \\
    \midrule
    
    \multicolumn{1}{l|}{Electric truck, include:} &  & & \\
    \multicolumn{1}{l|}{$\boldsymbol{\cdot}$ Base model} & \$250,000/vehicle & 10 yrs & $c^{veh}$: \$36,988/vehicle \\
    \multicolumn{1}{l|}{$\boldsymbol{\cdot}$ Battery} & \$150/kWh & 10 yrs & $c^{btr}$: \$22/kWh \\
    \midrule
    
    \multicolumn{1}{l|}{Charging infrastructure, include:} &  & & \\
    \multicolumn{1}{l|}{$\boldsymbol{\cdot}$ Charging station construction} & \$1,000,000/station & 20 yrs & $c^{ctr}$: \$106,781/station \\
    \multicolumn{1}{l|}{$\boldsymbol{\cdot}$ Power equipment} & \$200/kW & 20 yrs & $c^{cap}$: \$20/kW\\
    \multicolumn{1}{l|}{$\boldsymbol{\cdot}$ Chargers  (linearized)} & \$587/kW & 10 yrs & $c^{chg}$: \$87/kW\\
    \midrule

    \multicolumn{1}{l|}{Power delivery system, include:} &  & & \\
    \multicolumn{1}{l|}{$\boldsymbol{\cdot}$ Distribution lines} & \$1,200,000/mile & 30 yrs & $c^{lne}$: \$115,723/mile\\
    \multicolumn{1}{l|}{$\boldsymbol{\cdot}$ Substation upgrade (Fixed)} & \$4,600,000/upgrade & 25 yrs & $c^{upg,std}$: \$460,703/upgrade\\
    \multicolumn{1}{l|}{$\boldsymbol{\cdot}$ Substation upgrade (Variable)} & \$200,000/MW & 25 yrs & $c^{upg,var}$: \$20,031/MW\\
    \midrule
    
    \multicolumn{1}{l|}{Electricity (time-of-use)} & \multicolumn{3}{l}{$c^{ele}_t$: On-Peak (4PM-9PM): \$0.232/kWh} \\
    & \multicolumn{3}{l}{\hspace{6.3 mm} Mid-Peak (2PM-4PM \& 9PM-11PM): \$0.177/kWh} \\
    & \multicolumn{3}{l}{\hspace{6.3 mm} Off-Peak (all other hours): \$0.130/kWh} \\

    \bottomrule
\end{tabular}
\end{threeparttable}
\end{center}
\label{tab:cost}
\end{table}

This section discusses the details of the cost estimates, including initial investments, asset lifespans, and amortized costs, annualized at a 10\% interest rate. The cost of the electric truck system consists of the basic vehicle model \cite{alvarez2023tesla} and the battery \cite{iea2023trends}, with the base model featuring a 900kWh battery as per \cite{kothari2023tesla}. We assume a 10-year lifespan for both electric trucks and their batteries. Charging infrastructure costs include a \$1 million construction cost per site, assuming a 20-year lifespan. This estimate can be adjusted per station for more detailed models. The initial investment for power equipment is set at \$200 per kW, including the cost of transformers and other essential components like wires, conduits, and meters \cite{nelder2019reducing}, with an estimated lifespan of 20 years \cite{biccen2014lifetime}. The charger costs, based on data from \cite{crc2023access}, show a stepwise increase relative to the charger's power specification. For simplicity, we linearize these costs to a rate of \$587 per kW. Chargers are assumed to have a lifespan of 10 years. For the power delivery system, utilizing data from \cite{miso2022cost} and insights from a local utility company, we project the cost of distribution lines at approximately \$1.20 million per mile considering it is influenced by various factors like land needs and line types. For a substation upgrade, which includes a standard 28MVA transformer and necessary capacitor banks, the estimated cost is \$4.60 million per upgrade, again based on \cite{miso2022cost} and further consultations with the local utility. The variable cost of substation upgrades is estimated to be \$200/kW or \$200,000/MW, aligning with the cost of power equipment for charging stations. Electricity costs are based on time-of-use commercial rates adopted from \cite{ye2022decarbonizing}. Table \ref{tab:cost} summarizes these costs, lifespans, and their amortized equivalents.

\color{black}

\section{Results and Discussion}
\label{sec:results}
We solve the optimization model outlined in Section \ref{sec:method} under two different modes: 1) Hosting capacity examination and 2) Cost-effective regulation compliance, utilizing the data specified in Section \ref{sec:data}. The optimization model is formulated as a mixed-integer linear programming (MILP) problem and is solved with commercial solver Gurobi on a 32-core AMD machine.

\subsection{Hosting Capacity Examination}
\label{sec:hosting-capacity-examine}


Our first objective is to examine the electric drayage truck integration limits of the existing regional electric grid, specifically assessing how many electric drayage trucks the existing grid's remaining load hosting capacity can support. This involves solving the model under objective Mode 1: Hosting Capacity Examination, as detailed in Section \ref{sec:method-obj-mode1}, with the objective function of maximizing the number of electric drayage trucks. For this analysis, we temporarily do not consider substation upgrades to fully evaluate the existing grid's hosting capacity. However, it is also important to note that the grid's hosting capacity cannot be exclusively dedicated to electric drayage trucks due to competing demands from other sectors like residential, industrial, and other transportation modes. To approximate these external influences, we limit the available grid capacity in our analysis. Specifically, we look at scenarios where 20\%, 50\%, and 100\% of the grid's current capacity is allocated to electric drayage trucks. For the first two scenarios, we adjust the current remaining load hosting capacity of each substation, $P^{sbs}_k, \forall k \in K$, by multiplying 20\% and 50\%, respectively.

\begin{figure}[!h]
\centering
\captionsetup{justification=centering}
\includegraphics[width=0.90\textwidth, clip=true, trim = 0mm 0mm 0mm 0mm]{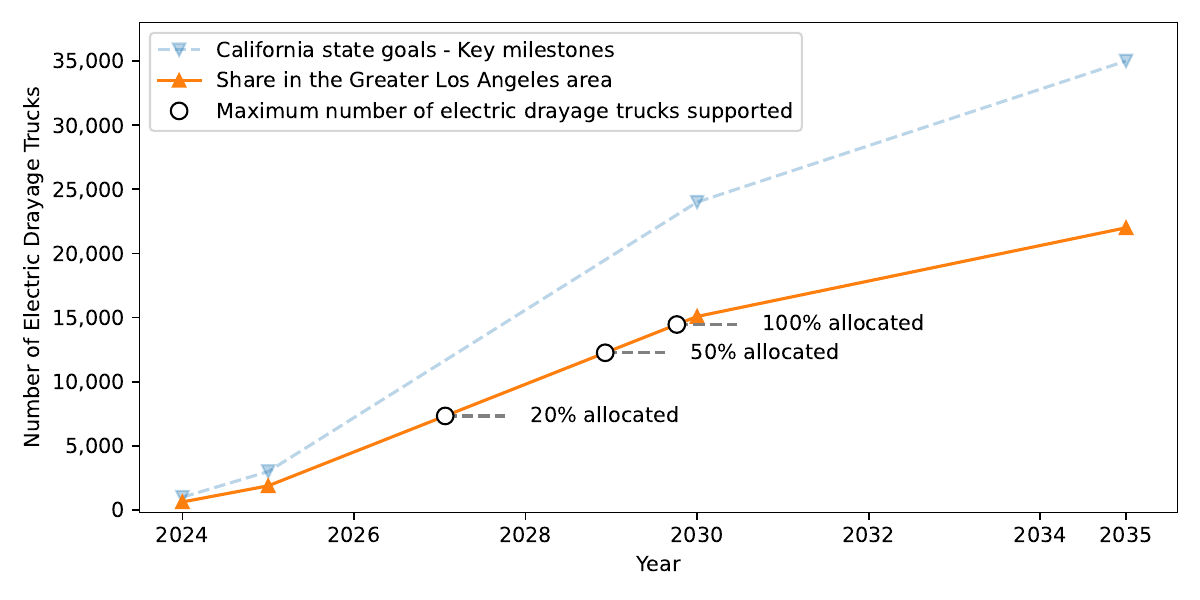}
\caption{Projected maximum number of electric drayage trucks supported, with $x$\% of the grid's current remaining load hosting capacity allocated to this sector ($x=$20, 50, and 100).}
\label{fig:results-capacity-examine}
\end{figure}

Figure \ref{fig:results-capacity-examine} illustrates the number of electric drayage trucks that can be supported under these scenarios, juxtaposed with California’s regulatory goals and the specific targets for the Greater Los Angeles area. Allocating 100\% of the remaining grid capacity to the drayage truck sector enables us to meet the 2030 goal, while a 50\% allocation might reach the goal around 2029. However, the possibility of allocating such high proportions of capacity exclusively to the drayage truck sector, given increasing demands from other sectors, is quite low. A more conservative 20\% allocation scenario indicates the goal could be met by 2027 without upgrades. However, considering that substation upgrades can take 5 to 15 years for planning, permitting, and completion, we could be already behind schedule in achieving the California's regulatory goals in electrifying the drayage trucks. Therefore, immediate and well-planned upgrades to the power delivery system in the greater Los Angeles area are crucial and should commence promptly.

\begin{figure}[!h]
\centering
\captionsetup{justification=centering}
\includegraphics[width=0.92\textwidth, clip=true, trim = 0mm 0mm 0mm 0mm]{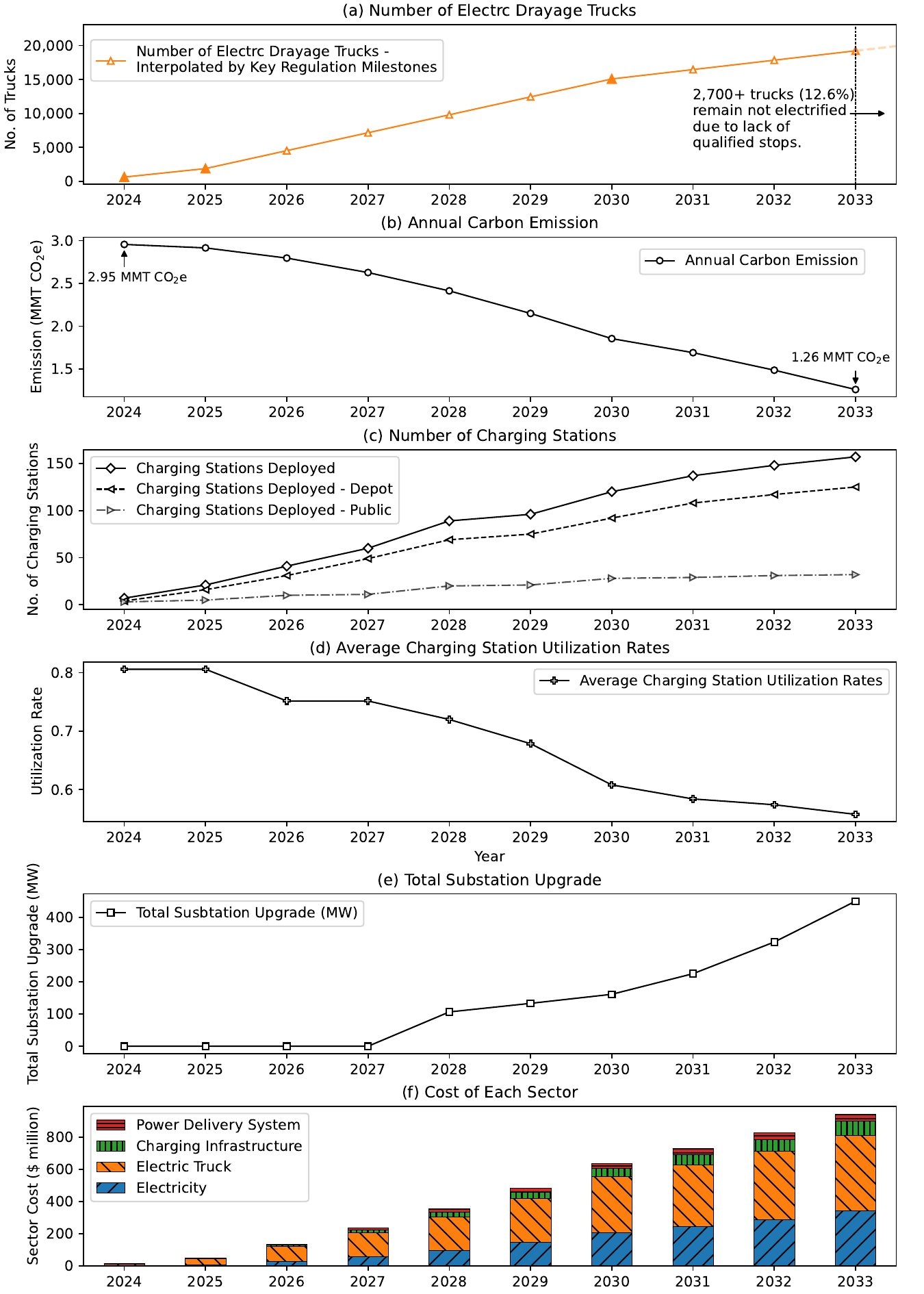}
\caption{Yearly development progress of the system under the most cost-effective scenario.}
\label{fig:progress}
\end{figure}

\subsection{Cost-Effective Regulation Compliance}
Amidst the pressing need for upgrades to the power delivery system, a critical question emerges: What is the most effective integrated strategy for developing multiple sectors in a coordinated manner? To answer this question, we focus on objective Mode 2: Cost-Effective Regulation Compliance in this subsection. Our approach begins optimistically by allowing a 100\% allocation of the current grid hosting capacity to the drayage truck sector. This setup helps us determine the timing for necessary upgrades to achieve the most cost-efficient planning. It's noteworthy that with any actual allocation below 100\%, the need for earlier upgrades becomes more pronounced. Additionally, to monitor the system's evolution on a yearly basis, we interpolate the Greater Los Angeles area's share of statewide milestones, breaking down the targets year by year as the $N_y$ variable in \eqref{eq:cons-truck-add} from Section \ref{sec:method-obj-mode2}. We also use electricity carbon intensity data $g_t$ from \cite{caiso2021eci} for our chosen representative day, which details total grid GHG emissions and power demands at 5-minute intervals.

\subsubsection{Development Progress}

In Figure \ref{fig:progress}, we track the annual development progress across various aspects of the system. Figure \ref{fig:progress} (a) presents the yearly targets for the Greater Los Angeles area, where hollow points are interpolated from the solid point milestones set by state regulations. Under the assumption that trucks' operational patterns remain consistent post-conversion, we observe that approximately 12.6\% of trucks lack adequate qualified stops for charging during their daily operations, likely due to intensive utilization and multiple shifts. These trucks may require alternative zero-emission solutions, such as fuel-cell technology, wireless charging infrastructure, or operational pattern adjustments. Our analysis primarily focuses on the drayage trucks scheduled for electrification up to and including the year 2033.

Figure \ref{fig:progress} (b) illustrates the reduction in GHG emissions in conjunction with the rising adoption of electric drayage trucks, taking into account the emissions from power generation. The GHG emissions, quantified in million metric tonnes (MMT) of CO$_2$, show a significant drop from 2.95 MMT in 2024 to 1.26 MMT by 2033, marking a reduction of over 57\%. It's important to note that our optimization model, as per objective function \eqref{eq:obj-2} in Section \ref{sec:method-obj-mode2}, does not directly use GHG emissions as a parameter to refine charging schedules. Therefore, there is potential for even greater reductions in GHG emissions through optimized charging schedules that actively consider real-time or prediction GHG emission \cite{wang2023predict}, coupled with a shift in the power grid from fossil fuel-based power plants to cleaner energy sources.

Figures \ref{fig:progress} (c) and (d) present the growth in the number of charging stations and their respective utilization rates. The utilization rate for charging station $j$ is defined as $UR = \left( \sum_{t \in T}P_{jt} \cdot \Delta t \right) / \left(P^{chs} \cdot T \right)$, representing the ratio of actual energy delivered to the trucks to the maximum energy that could be charged if the station operated at full capacity continuously. Initially, utilization rates are high but gradually reduce with the addition of more charging stations. Figure \ref{fig:progress} (e) outlines the optimal timeline for substation upgrades. In this scenario, where 100\% of the current grid capacity is allocated to electric drayage trucks, upgrades exceeding 100 MW are projected to start from 2028. It is crucial to note, however, that due to the lengthy planning and construction process of substation upgrades and the overly optimistic assumption of allocating 100\% of hosting capacity to the drayage truck sector, earlier initiation of substation upgrades should be strongly encouraged. Lastly, Figure \ref{fig:progress} (f) illustrates the cost distribution across different sectors. Here, the investment in electric trucks is the most significant, followed by the cost of electricity. This cost allocation underscores the financial expenses of transitioning to electric drayage trucks and the associated infrastructure needs.

\subsubsection{Geographical Expansion}

\begin{figure}[!h]
\centering
\captionsetup{justification=centering}
\includegraphics[width=0.97\textwidth, clip=true, trim = 0mm 0mm 0mm 0mm]{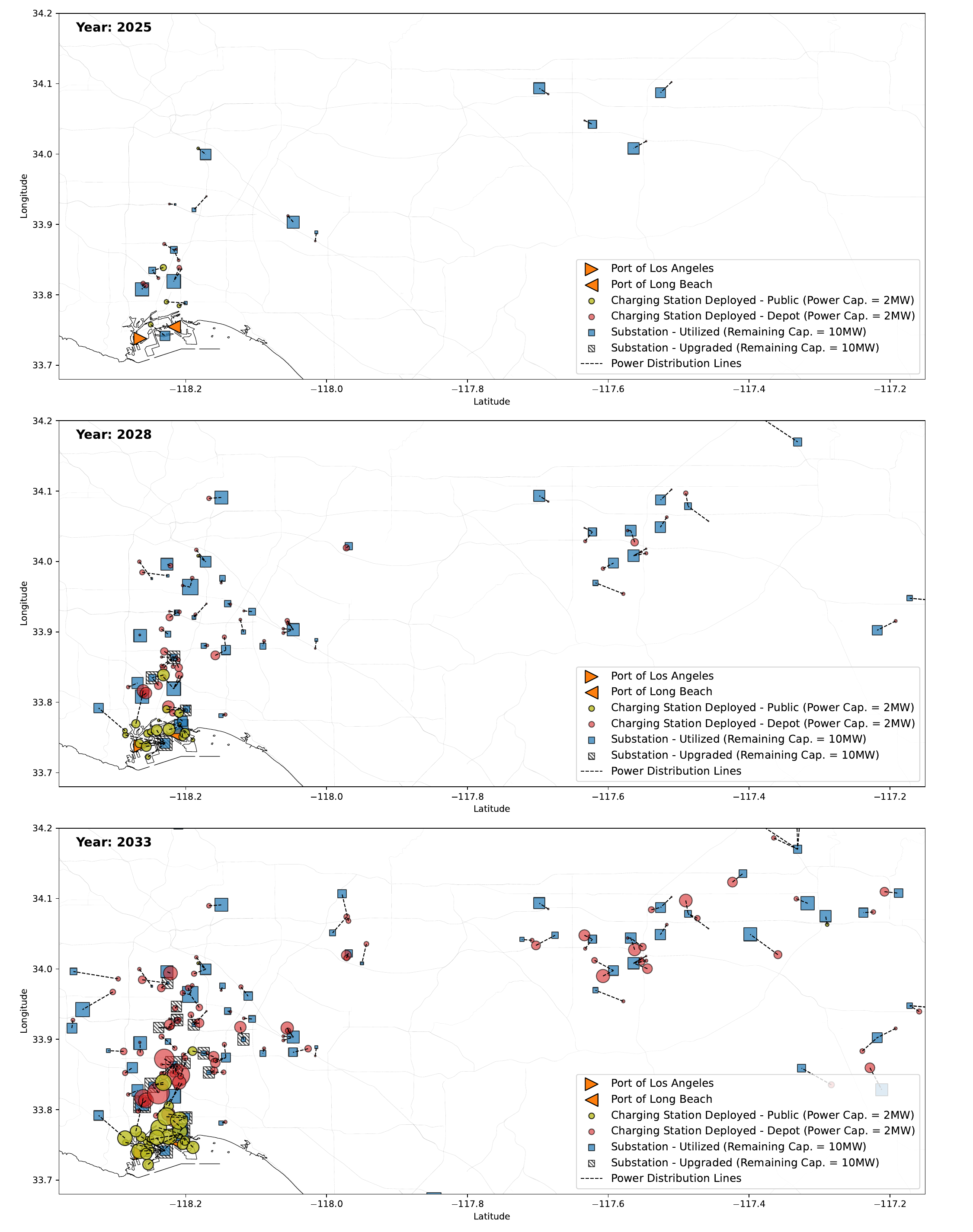}
\caption{Geographical distribution and expansion of charging stations, along with substation utilization and upgrades.}
\label{fig:results-geo-expansion}
\end{figure}

Beyond the numerical results, the spatial placement of charging stations, their corresponding substations, and the identification of substations requiring upgrades are key results. Figure \ref{fig:results-geo-expansion} illustrates the geographical expansion of these elements from 2025 to 2033, with a mid-point view in 2028. Notably, there is a significant increase in both the number and power capacity of charging stations near the port area, predominantly consisting of public charging stations. This area also experiences a high density of substation upgrades, indicative of substantial demand. Conversely, as we move inland, depots emerge as the most cost-effective charging locations due to longer stopping durations. Assuming the possibility of allocating 100\% of the current grid hosting capacity to the drayage truck sector, the model suggests that the inland grid capacity is adequate, with no further upgrades deemed necessary.

\subsubsection{Load Profiles}

\begin{figure}[!h]
\centering
\captionsetup{justification=centering}
\includegraphics[width=0.98\textwidth, clip=true, trim = 0mm 0mm 0mm 0mm]{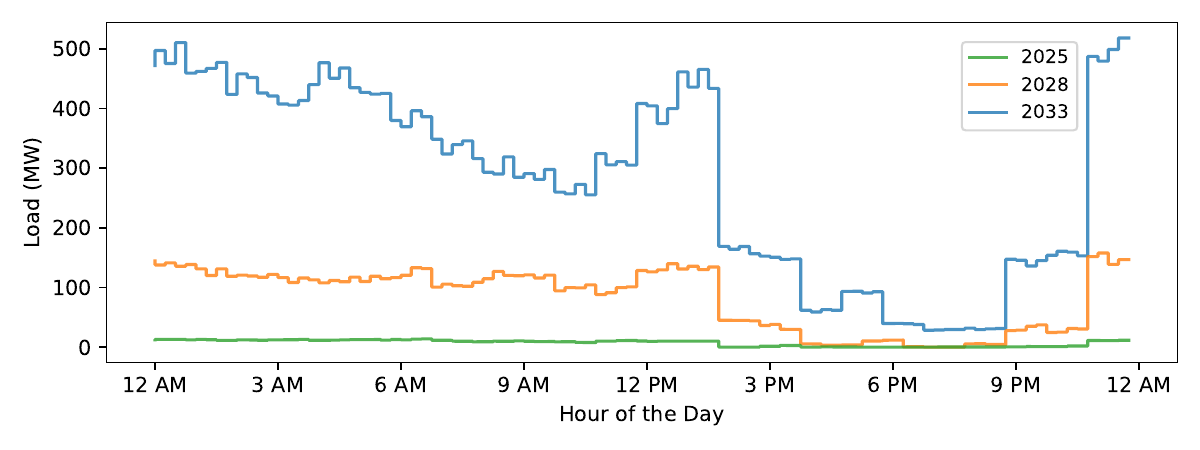}
\caption{Aggregated charging load profiles for drayage trucks in 2025, 2028, and 2033.}
\label{fig:results-load-profiles}
\end{figure}

Figure \ref{fig:results-load-profiles} illustrates the aggregated charging load profiles of the drayage truck sector for the years 2025, 2028, and 2033. Our model primarily optimizes charging schedules based on time-of-use electricity pricing, which is reflected in the load profiles that show charging activities strategically scheduled to avoid mid- and on-peak hours. Additionally, the model is also designed to minimize peak demand to mitigate costs associated with power system upgrades and charging station capacities. In this optimized scenario, the peak aggregated charging power exceeds 500 MW by 2033. However, considering potential variability in truck behavior, the actual peak demand might be even higher.

\subsubsection{Battery Capacity Requirements}

\begin{figure}[!h]
\centering
\captionsetup{justification=centering}
\includegraphics[width=0.98\textwidth, clip=true, trim = 0mm 0mm 0mm 0mm]{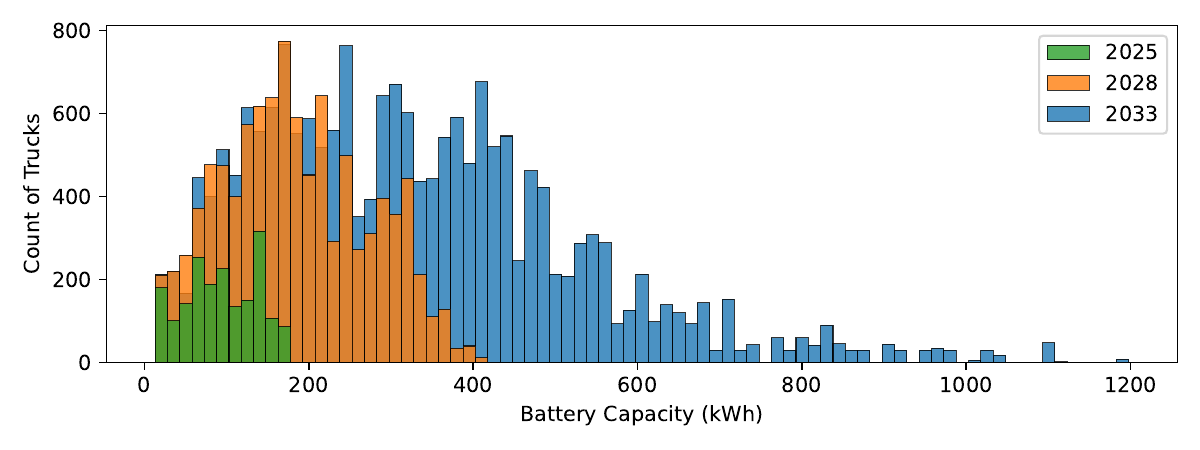}
\caption{Histograms of required battery capacities for electric drayage trucks in 2025, 2028, and 2033.}
\label{fig:results-battery-hist}
\end{figure}

For fleet owners and OEMs, understanding the necessary battery capacities for sustaining daily operations is crucial. Figure \ref{fig:results-battery-hist} shows the battery capacity requirements for the electric drayage truck fleet over the selected years. In the early stages of the transition, the battery capacity needs are relatively modest, indicating a prioritization of trucks with shorter daily routes and/or more frequent access to charging facilities. As the transition progresses, the average required battery capacity increases, signifying the need for larger batteries to support the electrification of additional drayage trucks. By 2033, the average battery capacity needed is projected to be around 328 kWh, with some trucks potentially requiring up to 1,200 kWh battery, surpassing the capacities of current leading models \cite{kothari2023tesla}.

\section{Conclusion}
\label{sec:conclusion}

This paper introduces a detailed optimization model designed to facilitate the transition to zero-emission drayage trucks, focusing on both infrastructure and operational requirements. Conducted in the context of the Greater Los Angeles area, the study underscores the critical need for extensive power system upgrades to align with California's ambitious zero-emission goals. The model provides insights into the most cost-effective and integrated strategies for developing the power delivery system, charging infrastructure, and electric drayage trucks, ensuring compliance with regulatory goals in a coordinated manner. It is noteworthy that the model's generic design allows for easy adaptation to other areas, simply by incorporating relevant regional data.

Considering the focus of this study on drayage trucks, future research could expand to explore the transition to zero-emission vehicles in other truck categories, including long-haul heavy-duty and local medium-duty trucks. Additionally, the concept of aggregating electric trucks as a virtual power plant, particularly with vehicle-to-grid capabilities, presents an interesting direction for future research. This approach has the potential to transform the role of electric trucks, enabling them not only to reduce emissions but also to enhance the stability and efficiency of the power grid, opening new frontiers in sustainable transportation and energy management.

\section*{Acknowledgments}
This study was made possible through funding received from University of California Office of the President and the University of California Institute of Transportation Studies from the State of California through the Public Transportation Account and the Road Repair and Accountability Act of 2017 (Senate Bill 1). The contents of this report reflect the views of the author(s), who is/are responsible for the facts and the accuracy of the information presented. This document is disseminated under the sponsorship of the State of California in the interest of information exchange and does not necessarily reflect the official views or policies of the State of California.

\bibliographystyle{apalike}  
\bibliography{references}  

\appendices

\section{GPS Data Processing}
\label{ap:gps}
We follow the procedures below to process the raw GPS data, transforming it into trajectories with uniform 15-minute intervals.

\begin{itemize}
    \item Step 1: Align GPS points with the road network using ArcGIS's ``Snap'' feature. This step corrects points that fall outside the road network due to measurement errors. We then use ArcGIS's network analyst toolset \cite{esri2022network} to find the shortest route between two consecutive raw GPS points, creating extra trajectory points that reflect the precise route.
    
    \item Step 2: Classify trajectory points into ``qualified stop'' or ``trip'' modes. A truck is considered in ``qualified stop'' mode if its speed is below 0.1 mph for 30 minutes or longer. Trucks in this mode are deemed eligible to use nearby charging stations.
    
    \item Step 3: Downsample the processed trajectories to create reduced-frequency trajectories with a uniform 15-minute interval. These downsampled points include travel distance to the next point and the proportion of time in ``qualified stop'' mode ($r^{stp}_{it}$), derived from the original high-frequency trajectory data. Figure \ref{fig:trajectories} illustrates the raw GPS data alongside the downsampled trajectories, with depot locations adjusted for confidentiality. The travel distance to the next point is used to calculate electricity consumption ($\Delta e^{con}_{it}$), assuming a consumption rate of 2 kWh per mile \cite{tesla2023semi}.
    
\end{itemize}

\end{document}